\documentclass[10pt,draft,reqno]{amsart}
 
\newtheorem{theorem}{Theorem}[section]

\newtheorem{proposition}[theorem]{Proposition}

\theoremstyle{definition}
\newtheorem{definition}[theorem]{Definition}
\newtheorem{example}[theorem]{Example}
\theoremstyle{remark}
\newtheorem{remark}[theorem]{Remark}
\numberwithin{equation}{section}
\newcommand{\half}{\frac{1}{2}}
\newcommand{\xx}{\mathbf{X}}
\newcommand{\ww}{\mathbf{W}}
\newcommand{\rr}{\mathbf{R}}

\begin{document}

\title[Super-exponents]{Brownian Super-exponents}

\author{Victor Goodman}
\address{Victor Goodman: Department of Mathematics, Indiana University, Bloomington, Indiana, U.S.A.}
\email{goodmanv@indiana.edu}



\begin{abstract}
We introduce a transform on the class of stochastic exponentials for
{\em d}-dimensional Brownian motions.  Each stochastic exponential
generates another stochastic exponential under the transform. The
new exponential process is often merely a supermartingale even in
cases where the original process is a martingale. We determine a
necessary and sufficient condition for the transform to be a
martingale process. The condition links expected values of the
transformed stochastic exponential to the distribution function of
certain time-integrals.
\end{abstract}

\maketitle

\section{Introduction}

If $\xx (t)$ is a {\em d}-dimensional progressively measurable
process and $\ww$ is a Brownian motion under a measure $P$, the {\em
stochastic exponential} determined by $\xx$ is the process

 \vspace{2ex}

$$ Z_{\mathbf{X}}(t) = \exp\left\{ \int_0^t \mathbf{X}(u)\cdot  d\mathbf{W}(u) -
\half \int_0^t ||\mathbf{X}(u)||^2du \right\}.$$
\vspace{4 ex}

\noindent The problem of checking whether $Z_{\xx}(t)$ is a true
martingale is important for the use of Girsanov's theorem.  Two
well-known sufficient conditions are due to Novikov and to Kazamaki;
see for example Revuz and Yor \cite{RY}.   Examples where the
process $Z_{\xx}(t)$ is strictly a supermartingale appear in Goodman
and Kim \cite{GK}, Levental and Skorohod \cite{LS}, and Wong and
Heyde \cite{WH}.

In their recent paper, Wong and Heyde \cite{WH} present a necessary
and sufficient condition for any stochastic exponential to form a
martingale process.  Their condition is formulated in terms of an
explosion time.  We consider a class of stochastic exponentials for
which their condition becomes more explicit. We begin with any
stochastic exponential and we describe a modification, or transform,
of it which generates another stochastic exponential.

 The transform involves a {\em time-integral} of the form

 $$\int_0^t ||\xx(u)||^2Z_{\xx}(u)du.$$

\noindent We derive a necessary and sufficient condition for the
transform to be a martingale.  Our condition is formulated in terms
of the distribution of time integrals, and we use the relation to
obtain bounds on the tail behavior of these distributions.

\begin{definition}\label{super}
Suppose that $\xx(t)$ is a progressively measurable process such
that for some $T>0$,

\begin{equation}\label{l2}
 P\left\{ \int_0^T ||\xx(u)||^2du < \infty  \right\} = 1
 \end{equation}

  If $Z_{\xx}(t)$ is the stochastic exponential generated by $\xx(t)$
   and $y >0$,  the associated
   \emph{super-exponent process} $Y_{\xx}(t)$ , defined for $t \le T$, is

\vspace{4 ex}

\begin{equation}\label{superdef}
 Y_{\xx ,y}(t)  = \frac{Z_{\xx}(t)}{y^{-1}+\half\int_0^t
||\xx(u)||^2Z_{\xx}(u)du}
\end{equation}
\vspace{4 ex}

\end{definition}
Notice from Equation (\ref{superdef}) that $ Y_{\xx ,y}(0) = y$. In
addition, $ Y_{\xx ,y}(t)$ is positive so that the random variable

$$\exp( Y_{\xx ,y}(t) )$$

\vspace{2 ex}

 \noindent is greater than one.  We show that this random variable
 has a finite expected
value which is less than or equal to $e^y$. This result is
surprising since $ Y_{\xx ,y}(t)$ is used as an exponent here.
According to Definition \ref{super}, $ Y_{\xx ,y}(t)$ itself
contains an exponential factor $Z_{\xx}(t)$.  For this reason, we
say that the process $ Y_{\xx ,y}(t)$ is a {\em Brownian
super-exponent}.
\vspace{3ex}

\section{Transform Properties}
\vspace{3ex}
\begin{proposition}\label{supersde}
Suppose that a progressively measurable process $\xx$ satisfies
condition (\ref{l2}).  Let $Y_{\xx,y}(t)$ denote the super-exponent
process in Definition \ref{super}.  Then for each $t \le T$,

\begin{equation}\label{ysde}
Y_{\xx,y}(t)=y + \int_0^tY_{\xx,y}(u)\xx(u)\cdot  d\ww(u) - \half
\int_0^t||\xx(u)||^2Y_{\xx,y}^2(u)du
\end{equation}
\vspace{2ex}

\noindent Moreover, the process
\begin{equation}
\exp\left(Y_{\xx,y}(t)\right)
\end{equation}
\vspace{2ex}

\noindent is a positive supermartingale on the interval $0 \le t \le
T$.

\noindent In addition, the process \vspace{3ex}
\begin{equation}\label{newse}
\tilde Z(t) = \exp\left(Y_{\xx,y}(t) - y \right)
\end{equation}

\vspace{3ex}

 \noindent is a stochastic exponential for $\ww$.  This
stochastic exponential is generated by the d-dimensional process
\begin{equation}
Y_{\xx,y}(t)\xx(t)
\end{equation}

\end{proposition}
\vspace{4ex}

\begin{proof} It follows from the definition of $Z_{\xx}(t)$ that

$$dZ_{\xx} = Z_{\xx}\xx\cdot d\ww \mbox{\hskip.5in and\hskip.5in}
d\int_0^t ||\xx||^2Z_{\xx}du= ||\xx||^2Z_{\xx}dt$$

\noindent Direct calculation shows that

\begin{equation}
\begin{aligned}
dY_{\xx,y}=\frac{dZ_{\xx}}{y^{-1}+\half\int_0^t ||\xx||^2Z_{\xx}du}
&+Z_{\xx}\ d(y^{-1}+\half\int_0^t ||\xx||^2Z_{\xx}du )^{-1}\\
&=Y_{\xx,y}\xx\cdot \ww -\half ||\xx||^2 Y_{\xx,y}^2dt\\
\end{aligned}
\end{equation}

From this equation we see that $Y_{\xx,y} - y$ is the sum of the
It\^ o integral of $Y_{\xx,y}\xx$ and the elementary integral of
$-\half ||Y_{\xx,y}\xx||^2 $.  This establishes Equation
(\ref{ysde}). It follows immediately from Equation (\ref{ysde}) that
$Y_{\xx,y}-y$ is the exponent of a stochastic exponential.
Therefore, the process
$$\exp\left(Y_{\xx,y}(t)-y \right)$$
is a positive local martingale.  It is well known that a positive
local martingale is a supermartingale; see, for instance, Karatzas
and Shreve \cite{KS}.  In addition, Equation (\ref{newse}) is a
direct consequence of Equation (\ref{ysde}) and the definition of
stochastic exponential processes.

\end{proof}

\begin{theorem}\label{dist}  Suppose that  $\xx(t)$ is a deterministic
function such that for some $T > 0$
$$\int_0^T ||\xx(u)||^2du < \infty.$$

\noindent Let $Z_{\xx}(t)$ and $Y_{\xx,y}$ denote the stochastic
exponential and super-exponent process generated by $\xx(t)$.  Then
for each non-negative measurable function $G(u)$, $u > 0$, and $t
<T$,

 \vspace{2ex}

\begin{equation}\label{timedist}
\begin{aligned}
&\phantom{xxxxxxxxx}  E[G(Y_{\xx,y}(t)) \exp (Y(t) - y)]\\
 &= E\left[\  G(\frac{ Z_{\xx}(t)}{y^{-1} -
\half\int_0^t||\xx(u)||^2Z_{\xx}(u)du}) \ ;\
\int_0^t||\xx(u)||^2Z_{\xx}(u)du
< \frac{2}{y}\  \right]\\
\end{aligned}
 \end{equation}
\vspace{2ex}

\end{theorem}

\begin{proof}  For $N= 1,2,\dots$ let $\tau_N$ be the stopping time defined by

$$\tau_N = \inf \{t \le T\ :\ Y_{\xx,y}(t) \ge N \}$$
\vspace{2ex}
 \noindent It follows from Equation (\ref{ysde}) that
$$Y_{\xx,y}(t\wedge \tau_N) - y  =  \int_0^{t\wedge
\tau_N}Y_{\xx,y}(u)\xx(u)\cdot d\ww(u) - \half \int_0^{t\wedge
\tau_N}||\xx(u)||^2Y_{\xx,y}^2(u)du$$

\noindent From this equation we see that $\exp(Y_{\xx,y}(t\wedge
\tau_N) - y)$ is another stochastic exponential which is generated
by

$$Y_{\xx,y}(u)1_{\{u < \tau_N\}}\xx(u) $$

Since this process is uniformly bounded in $L^2[0,T]$, it satisfies
Novikov's condition.  It is well known (see Karatzas and Shreve
\cite{KS}) that the associated stochastic exponential is a
martingale. We apply Girsanov's Theorem to change measure using the
Radon-Nykodym derivative

$$\Lambda (T) = \exp(Y_{\xx,y}(T\wedge \tau_N) - y)$$
The probability measure $Q_N$ is given by
$$\frac{dQ_N}{dP} = \Lambda (T) $$
Then with respect to $Q_N$ the process
$$\mathbf{ \tilde W}(t) = \ww(t) - \int_0^{t\wedge \tau_N} Y_{\xx,y}(u)\xx(u)du$$
is a Brownian motion for $t \le T$.  Since $Y_{\xx,y}(t)$ is a
strong solution to equation (\ref{ysde}), we may consider its SDE
with respect to the Brownian motion $\mathbf{ \tilde W}$:

For $t < \tau_N$

\begin{equation}\label{tildesde}
\begin{aligned}
dY_{\xx,y} = &Y_{\xx,y}\xx \cdot d\ww  - \half||\xx||^2Y_{\xx,y}^2dt\\
= &Y_{\xx,y}\xx \cdot \{\ d\mathbf{\tilde W} +Y_{\xx,y}\xx dt  \}
 - \half||\xx||^2Y_{\xx,y}^2dt\\
= &Y_{\xx,y}\xx \cdot d\mathbf{\tilde W}  + \half||\xx||^2Y_{\xx,y}^2dt\\
\end{aligned}
\end{equation}

\noindent Now we have an explicit solution to the SDE in equation
(\ref{tildesde}):
\begin{equation}\label{tildeyyy}
Y_{\xx,y}(t)  = \frac{\tilde Z_{\xx}(t)}{y^{-1}-\half\int_0^t
||\xx(u)||^2\tilde Z_{\xx}(u)du}
\end{equation}
In this equation,  $\tilde Z_{\mathbf{X}}(t)$ denotes the stochastic
exponential (generated by $\mathbf{X}$) with respect to  the
Brownian motion $\mathbf{ \tilde W}$.  Now we consider
\begin{equation}\label{ident}
\begin{aligned}
 &E[G(Y_{\xx,y}(t)) \exp(Y_{\xx,y}(t) - y)1_{\{ t < \tau_N\}}]\\
 &= E[ G(Y_{\xx,y}(t))\Lambda (T)1_{\{ t < \tau_N\}}]\\
&= E_{Q_N}[G(Y_{\xx,y}(t))1_{\{ t < \tau_N\}} ]\\
&= E_{Q_N}[G(\frac{\tilde Z_{\xx}(t)}{y^{-1}-\half\int_0^t
||\xx(u)||^2\tilde Z_{\xx}(u)du})1_{\{ t < \tau_N\}} ]\\
\end{aligned}
\end{equation}
Here we used the identity for $Y_{\xx,y}$ in Equation
(\ref{tildeyyy}).

Moreover, from Equation (\ref{tildeyyy}) we also have

$$ t < \tau_N \mbox{\hskip.3in iff.\hskip.3in} \max_{s\le t}\frac{\tilde
Z_{\xx}(s)}{y^{-1}-\half\int_0^s ||\xx(u)||^2\tilde Z_{\xx}(u)du}) <
N$$
This allows us to write the last expected value in Equation
(\ref{ident}) as

$$E[G(\frac{ Z_{\xx}(t)}{y^{-1}-\half\int_0^t
||\xx(u)||^2\ Z_{\xx}(u)du})\ ;\ \max_{s\le t}\frac{
Z_{\xx}(s)}{y^{-1}-\half\int_0^s ||\xx(u)||^2 Z_{\xx}(u)du}) < N ]$$

\noindent since the integrand involves only the distribution of a
Brownian motion for each choice of $N$.  The limit of this expected
value as $N\to\infty$ is

$$E[G(\frac{ Z_{\xx}(t)}{y^{-1}-\half\int_0^t
||\xx(u)||^2\ Z_{\xx}(u)du})\ ;\ \half\int_0^t ||\xx(u)||^2
Z_{\xx}(u)du) < y^{-1} ]$$

Since the limit of the first expected value in  Equation
(\ref{ident}) is
$$E[G(Y_{\xx,y}(t)) \exp(Y_{\xx,y}(t) - y)]$$
the theorem is proved.
\end{proof}
\vspace{3ex}
\section{Examples Using the Transform}
\vspace{4ex}

\begin{proposition}  Suppose that  $\xx(t)$ is a deterministic
function such that
$$\int_0^t ||\xx(u)||^2du$$

\noindent is strictly increasing and finite for $t \le T < \infty$.
Let $Z_{\xx}(t)$ and $Y_{\xx,y}$ denote the stochastic exponential
and super-exponent process generated by $\xx(t)$.  Then the process
$$\exp( Y_{\xx,y}(t))$$

\noindent is a strict supermartingale for $t \le T$.  Moreover,

\begin{equation}\label{integral}
E[ \exp( Y_{\xx,y}(t)) ] = e^y Pr
\left\{\int_0^t||\xx(u)||^2Z_{\xx}(u)du < \frac{2}{y}\ \right\}
\end{equation}
\end{proposition}

\vspace{4ex}

\begin{proof}
We apply Theorem \ref{dist} using the choice $G(u)\equiv 1$.
Equation (\ref{timedist}) becomes
$$ E[\exp (Y_{\xx,y}(t) - y)] = Pr\left\{\int_0^t||\xx(u)||^2Z_{\xx}(u)du
 < \frac{2}{y}\ \right\},$$
and Equation (\ref{integral}) follows.  Now since each $ Z_{\xx}(u)$
is a log normal random variable, the process
\begin{equation}\label{time}
\int_0^t||\xx(u)||^2Z_{\xx}(u)du
\end{equation}
has strictly increasing sample paths.  It follows that the right
hand expression in Equation (\ref{integral}) is strictly decreasing.
Therefore, $\exp (Y_{\xx,y}(t))$ is a strict supermartingale.

\end{proof}

\vspace{4ex}

\begin{remark}  Equation (\ref{integral}) provides a useful tool for
investigating the distribution of a time integral given by Equation
(\ref{time}). Since each super-exponent

$$Y_{\xx ,y}(t)  = \frac{Z_{\xx}(t)}{y^{-1}+\half\int_0^t
||\xx(u)||^2Z_{\xx}(u)du}$$

\noindent is point-wise increasing as a function of $y$, it follows
from the identity
$$Pr\left\{\int_0^t||\xx(u)||^2Z_{\xx}(u)du
 < a\ \right\} = \exp (-\frac{2}{a})E[\exp(Y_{\xx ,2/a}(t)   )]$$
that the distribution function is the product of a decreasing
function of $a$ and the explicit factor $\exp (-2/a)$ .

It is not known whether $\exp(Y_{\xx ,\infty}(t))$ has finite
expectation. A finite expected value would produce sharp estimates
for the lower tail probability of (\ref{time}).  We conjecture that

$$E[\ \frac{2Z_{\xx}(t)}{\int_0^t
||\xx(u)||^2Z_{\xx}(u)du}\ ] = \infty.$$
\end{remark}

\vspace{3ex}
\begin{example} In the case of $d=1$ the choice $X(t)\equiv
\sigma$ specializes the time integral in (\ref{time}) to a time
integral of {\em geometric Brownian motion:}

\begin{equation}\label{fin}
\int_0^t \exp(\sigma W(u) - \sigma^2u/2)du
\end{equation}

Expected values involving related time integrals appear in
computational problems of financial mathematics. Consequently,
distribution properties of these time integrals have been studied by
many authors; see Dufresne \cite{D}, Geman and Yor \cite{GY}, Rogers
and Shi \cite{RS}, and Goodman and Kim \cite{GK}.

Although most works have used analytic techniques to express the
distribution in various integral forms, in Goodman and Kim \cite{GK}
martingales techniques are used exclusively.   A special case of Equation (\ref{integral})
  appears in  \cite{GK}, Theorem 4.1:
$$Pr\left\{\int_0^t\exp(W(u)-u/2)du
 \le a\ \right\}$$
 $$ = \exp (-\frac{2}{a})E[\exp(\frac{2\exp(W(t)-t/2)}
 {a+\int_0^t\exp(W(u)-u/2)du}   )]$$
 \vspace{2ex}

\noindent The right hand expression for the distribution can be differentiated
with respect to $a$.  Consequently, it is shown in \cite{GK} that
the density function  multiplied by $a^2/2$ equals the
difference between two distribution functions of time integrals of slightly different geometric Brownian motions.
\end{example}
\vspace{3ex}
\begin{example}

\end{example}  In contrast to deterministic choices for $\xx(t)$,
where the stochastic exponential
$$\exp(Y_{\xx,y}(t))$$
is never a martingale, stochastic choices for $\xx$  may produce
martingales.  Of course, the introduction of a stopping time, as we have seen in the
proof of Theorem \ref{dist}, may produce a martingale.  In other cases, stopping times are
not required.

Consider the example of $X(t)= \cos(W(t))$, again in
the case $d=1$.  Then

$$Z_{\xx}(t) = \exp(\int_0^t \cos(W(u))dW(u) - \half\int_0^t \cos^2(W(u))du)$$
$$=\exp(\sin(W(t)) + \half\int_0^t[\sin(W(u) - \cos^2(W(u))]du)$$
is a bounded random variable.  Therefore, its super-exponent,
$Y_{\cos(W),y}(t)$ is also bounded.  Then since the local martingale
$$\exp(Y_{\cos(W),y}(t))$$
is also bounded, it is a martingale.  It is of interest then to know
when a super-exponent generates a martingale process.

 \vspace{4ex}
\section{The Martingale Condition}

 \vspace{4ex}

Theorem 1 of Wong and Heyde \cite{WH} identifies a necessary and sufficient condition for a progressively measurable process $\tilde \xx$ to generate a martingale stochastic exponential  process.  For completeness, we state their result here.

 \vspace{4ex}

 \begin{proposition} {\em (\cite{WH}, Proposition 1)} Consider a d-dimensional progressively measurable process $\tilde \xx(t) = \xi(\ww(\cdot) , t)$.  Then there will also exist a d-dimensional progressively measurable process

 $$\tilde \rr(t) = \xi(\ww(\cdot) +\int_0^{\cdot} \tilde \rr (u)du , t)$$

 \noindent defined possibly up to an explosion time $\tau^{M_\rr}$ where

 $$\tau^{M_R} = \inf \left\{ t \in \mathbb{R}^+ : M_\rr (t) = \int_0^t||\tilde \rr(u)||^2du = \infty \right\}$$

 \end{proposition}

  \vspace{4ex}

  \begin{theorem} {\em (\cite{WH}, Theorem 1)} Consider $\tilde \xx(t)$ and $\tilde \rr(t)$ as defined in Proposition 4.1.  The stochastic exponential $Z_{\tilde \xx}(T)$  satisfies

  $$P(\tau^{M_\rr} > T) = E_{P}[Z_{\tilde \xx}(T) ]$$
   \vspace{1ex}

  \noindent and hence is a martingale if and only if $P(\tau^{M_\rr} > T)=1$.

 \end{theorem}

  \vspace{4ex}

 We apply Theorem 1 of \cite{WH} using $\tilde \xx(t) = Y_{\xx,y}(t)\xx(t)$.   That is, our generating process is the one in Proposition \ref{supersde} where the stochastic exponential process is

 $$\exp(Y_{\xx,y}(t) - y ).$$

    \vspace{1ex}

\noindent  We first show that each generating process $\xx$ implicitly defines
 another process $\xx'$.  This allows us to identify the process $\tilde \rr(t)$.

 \vspace{4ex}


 \begin{proposition}\label{prime}
Suppose that a  d-dimensional progressively measurable process
$\xx(t)$ satisfies

$$Pr\left(\int_0^T ||\xx(u)||^2du < \infty\right) = 1$$
for some $T > 0$.  Then there exists another progressively
measurable process $\xx'(t)$,  so that  if $\tilde \xx(t) := Y_{\xx,y}(t)\xx(t)1_{\{t\le T\}}$ in Proposition 4.1, then the process $\tilde \rr(t)$ of the proposition satisfies

$$\tilde \rr(t) = \frac{ Z_{\xx'}(t)}{y^{-1}-\half\int_0^t ||\xx'(u)||^2 Z_{\xx'}(u)du}\xx'(t)$$
for all $t < \tau^{M_\rr}$.

\noindent Moreover,

 $$\tau^{M_R} = \inf \left\{ t \in \mathbb{R}^+ :  \int_0^{t\wedge T}||\xx'(u)||^2Z_{\xx'}(u)du = 2/y \right\}$$

\end{proposition}

\vspace{4ex}


\begin{proof} We follow the proof of Proposition 4.1.
 Let
 $$\tilde \xx(t) := Y_{\xx,y}(t)\xx(t)1_{\{t\le T\}}.$$
For each $N = 1,2, \dots$ we define a sequence of stopping times by

$$\tau_N = \inf\{t \in \mathbb{R}^+ : \int_0^t  Y_{\xx,y}^2(u)||\xx(u)||^21_{\{u\le T\}}du \ge N\}$$

It follows from Equation (\ref{ysde}) that
$$Z_{\tilde \xx}(t\wedge \tau_N) = \exp(Y_{\xx,y}(t\wedge \tau_N) - y)$$
 forms a martingale.  As in the proof of Theorem \ref{dist}, we apply
Girsanov's theorem using the Radon-Nikodym derivative
$$\Lambda (T) = \exp(Y_{\xx,y}(T\wedge \tau_N) - y)$$
to obtain the probability measure $Q_N$ where
$$dQ_N = \Lambda (T)dP. $$
With respect to the measure $Q_N$, the process
$$\mathbf{ \tilde W}(t) = \ww(t) - \int_0^{t\wedge \tau_N} Y_{\xx,y}(u)\xx(u)du$$
is a Brownian motion.  Hence, on the set $\{ t \le \tau_N \wedge T \}$ we have

$$\tilde \xx(t) = \xi(\mathbf{ \tilde W}(\cdot) +  \int_0^{\cdot} Y_{\xx,y}(u)\xx(u)du, t)$$
That is,

\begin{equation}\label{heyde}
Y_{\xx,y}(t)\xx(t) = \xi(\mathbf{ \tilde W}(\cdot) +  \int_0^{\cdot} Y_{\xx,y}(u)\xx(u)du, t)
\end{equation}

\vspace{4ex}

Now the process $Y_{\xx,y}(t)$ can also be described in terms of the Brownian motion $\mathbf{ \tilde W}$.
The calculations in Equation (\ref{tildesde}) also apply to the stochastic case.   Equation (\ref{tildeyyy}) gives an explicit  formula for $Y_{\xx,y}$:

\begin{equation}\label{form}
Y_{\xx,y}(t)  = \frac{\tilde Z_{\xx}(t)}{y^{-1}-\half\int_0^t
||\xx(u)||^2\tilde Z_{\xx}(u)du}
\end{equation}

We see that each term of Equation (\ref{heyde}) is a functional of $\xx$ and the Brownian motion $\tilde \ww$.
This  demonstrates the existence of a process  $\xx$ so that (\ref{heyde}) and (\ref{form}) hold up to a time $\tau_N$ defined by the integral of $Y_{\xx,y}(u)\xx(u)$ , using the Brownian motion $\tilde \ww$.

 Therefore, using the identical distribution of $\ww$ and the (original) measure $P$, we see that there exists a progressively measurable process  $\xx'(t)$ so that

\begin{equation*}
\frac{ Z_{\xx'}(t)}{y^{-1}-\half\int_0^t ||\xx'(u)||^2 Z_{\xx'}(u)du}\xx'(t)
= \xi(\mathbf{  W}(\cdot) +  \int_0^{\cdot} Y_{\xx',y}(u)\xx'(u)du, t)
\end{equation*}

\vspace{4ex} Here, we have abbreviated the complete expression on
the right hand side using (\ref{form}) to provide the notation. That
is, $Y_{\xx',y}$ denotes the expression in Equation (\ref{form}) but
{\em in the original Brownian motion} and $\xx$ is replaced by the
process $\xx'$.

 As $N\to\infty$ the stopping time  $\tau_N$ increases
to the stopping time

$$\tau = \inf\{t \le T : \int_0^t  Y_{\xx',y}^2(u)||\xx'(u)||^2du =\infty\}$$

\noindent By construction, the new process $\xx'$ satisfies

$$\int_0^T ||\xx'(u)||^2du < \infty \mbox{\hskip.3in a.\
s.\hskip.3in and\hskip.3in} \xx'(u) = 0 \mbox{\hskip.1in
for\hskip.1in} u>T$$

  Therefore, the process $Y_{\xx',y}$ (again, defined as in
  (\ref{form})) is bounded along each sample path up to the time where
  its denominator first hits zero. This defines the stopping time $\tau^{M_R}$ of the Proposition.
\end{proof}

\vspace{4ex}
\begin{theorem}  Suppose that $\xx(t)$ and $\xx'(t)$ are d-dimensional processes as defined in Proposition \ref{prime}.  Then the super-exponent process $Y_{\xx,y}(t)$ satisfies

\begin{equation}\label{mart}
E[ \exp(Y_{\xx,y}(t) - y)] = Pr\left\{\int_0^t ||\xx'(u)||^2Z_{\xx}(u)du < 2/y\right\}
\end{equation}
for $t  \le T$.
\end{theorem}

\vspace{4ex}

\begin{proof}  From Theorem 4.2 and Proposition \ref{prime} we have

$$ E[ \exp(Y_{\xx,y}(t) - y)] = Pr\left \{ \tau^{M_R} > t \right\}$$

$$=Pr\left \{  \int_0^t  Y_{\xx',y}^2(u)||\xx'(u)||^2du <\infty\right\}$$

$$=Pr\left \{  \int_0^t  \frac{ Z_{\xx'}(u)}{y^{-1}-\half\int_0^u ||\xx'(r)||^2 Z_{\xx'}(r)dr}||\xx'(u)||^2du <\infty\right\}$$

$$=Pr\left\{\int_0^t ||\xx'(u)||^2Z_{\xx}(u)du < 2/y\right\} $$

\end{proof}

 \vspace{4ex}

\bibliographystyle{amsplain}

\begin{thebibliography}{99}




\bibitem{D} Dufresne, D.:  The integral of geometric Brownian motion; {\em Adv. in Appl. Probab.}
\textbf{33} (2001) 223-241

\bibitem{GY} Geman, H, and Yor, M.: Asian Options, Bessel Processes and Perpetuities;
 {\em Math. Finance }\textbf{2} (1993) 349-375

\bibitem{GK}  Goodman, V. and Kim, K.:  Exponential martingales and
time integrals of Brownian motion  {\em preprint}

\bibitem{KS} Karatzas, I., and Shreve, S.: {\em Brownian Motion and Stochastic Calculus.}
 Springer-Verlag, New York, 1991

\bibitem{K} Kim, K.: Moment Generating function of the inverse of integral of geometric Brownian Motion;
{\em Proc. Amer. Math. Soc.} \textbf{132} (2004) 2753-2759

\bibitem{LS} Levental, S. and Skorohod, A. V.:
A necessary and sufficient condition for absence of arbitrage with
tame portfolios; {\em Ann. Appl. Prob.} \textbf{5}  (1995) 906-925

\bibitem{RY} Revuz, D. and Yor, M.: {\em Continuous Martingales and
Brownian Motion, 3rd edn.} Springer-Verlag, New York, 1999

\bibitem{RS} Rogers, L.C.G., and Shi, Z.:
The value of an Asian option {\em  J. Appl Appl. Probab. }
\textbf{32} (1995) 1077-1088

\bibitem{WH} Wong, B. and Heyde, C.C.:
On the martingale property of stochastic exponentials; {\em J.
Appl.Probab.} \textbf{41} (2004) 654-664

\bibitem{Y} Yor, M.: On some exponential functionals of Brownian
motion; {\em Adv. in Appl. Probab.} \textbf{24}  (1992) 509-531

\end{thebibliography}

\end{document}